\documentstyle[11pt,leqno]{article}
\newtheorem{th}{Theorem}[section]
\newtheorem{lem}[th]{Lemma}
\newtheorem{pro}[th]{Proposition}
\newtheorem{co}[th]{Corollary}

\title{\Large \bf {Geometry of Quaternionic K\"ahler connections with torsion }}
\author{{\sc Stefan Ivanov}
\thanks{The author is supported by Contract MM 809/1998 with the
Ministry of Science and Education of Bulgaria and by Contract 238/1998 with the
University of Sofia "St. Kl. Ohridski".}}
\date{}
\begin{document}
\maketitle
\thispagestyle{empty}
\vspace{2mm}
\begin{center}
{\sc Department of Mathematics\\ University of Sofia\\ "St. Kl. Ohridski" }
\end{center}
\vspace{5mm}

\begin{abstract}
The target space of a (4,0) supersymmetric two-dimensional sigma model with 
Wess-Zumino term has a connection with totally skew-symmetric torsion and  
holonomy contained in Sp(n).Sp(1), QKT-connection. We study  the 
geometry of QKT-connections. We find conditions to the existence of a 
QKT-connection and prove that if it exists it is unique. 
Studying conformal transformations we obtain a lot of (compact)
  examples
of QKT manifolds. We present a (local) description of 4-dimensional 
homogeneous QKT structures relying on the known result of naturally reductive 
homogeneous Riemannian manifolds. We consider Einstein-like QKT manifold and 
find  closed relations with Einstein-Weyl geometry in dimension four.
\\[8mm]
{\bf Running title:} Quaternionic K\"ahler with torsion
\\[8mm]
\\[8mm]
\end{abstract}

\section{Introduction and statement of the results}

An almost hyper complex structure on a 4n-dimensional manifold $M$ is a triple
$H=(J_{\alpha}), \alpha=1,2,3$, of almost complex structures
$J_{\alpha}:TM\rightarrow TM$ satisfying the quaternionic identities
$J_{\alpha}^2=-id$ and $J_1J_2=-J_2J_1=J_3$. When each $J_{\alpha}$ is a complex
structure, $H$ is said to be a hyper complex structure on $M$.

An almost quaternionic structure on $M$ is a rank-3 subbundle $Q \subset End(TM)$
which is locally spanned  by almost hypercomplex structure $H=(J_{\alpha})$; such a
locally defined triple $H$ will be called an admissible basis of $Q$. A linear
connection $\nabla$ on $TM$ is called quaternionic connection if $\nabla$ preserves
$Q$, i.e. $\nabla_X\sigma\in \Gamma (Q)$ for all vector fields $X$ and smooth
sections $\sigma \in \Gamma (Q)$. An almost quaternionic structure is said to be a
quaternionic if there is a torsion-free quaternionic connection. A $Q$-hermitian
metric is a Riemannian metric which is Hermitian with respect to each almost complex
structure in $Q$. An almost quaternionic (resp. quaternionic) manifold  with
Q-hermitian metric is called an almost quaternionic Hermitian (resp.
quaternionic hermitian) manifold

For $n=1$ an almost quaternionic structure is the same as an oriented conformal
structure and it turns out to be always quaternionic. When $n\ge 2$, the existence
of torsion-free quaternionic connection is a strong condition which is equivalent to
the 1-integrability of the associated GL(n,H)SP(1) structure \cite{Bon,Ob,Sal4}. 
If the Levi-Civita
connection of a quaternionic hermitian  manifold $(M,g,Q)$
is a quaternionic connection then $(M,g,Q)$ is called Quaternionic K\"ahler (briefly
QK). This condition is equivalent to the statement that the holonomy group of $g$ is
contained in SP(n).SP(1) \cite{A1,A2,S1,S2,Ish}.  If on a QK manifold there exist an
admissible basis $(H)$ such that each almost complex structure 
$(J_{\alpha})\in(H), \alpha =1,2,3$   is parallel with respect to
the Levi-Civita connection
 then the manifold is called hyper K\"ahler (briefly HK). In this
case the holonomy group of $g$ is contained in SP(n).

The notions of quaternionic manifolds arise in a natural way from the theory of
supersymmetric sigma models. The geometry of the target space of two-dimensional 
sigma models with extended supersymmetry is described by the properties of a 
metric connection with torsion \cite{HP1,HP2}. 
The geometry of (4,0) supersymmetric two-dimensional
sigma models without Wess-Zumino term (torsion) is a hyper K\"ahler manifold. In the
presence of torsion the geometry of the target space becomes hyper K\"ahler with
torsion (briefly HKT) \cite{HP3}. This means that the complex structures
$J_{\alpha},
\alpha=1,2,3$,  are parallel with respect to a metric quaternionic connection with
totally skew-symmetric torsion \cite{HP3}. Local (4,0)
supersymmetry requires that the target space of two dimensional sigma models with
Wess-Zumino term be either HKT or quaternionic K\"ahler with torsion (briefly QKT)
\cite{Hish} which means that the quaternionic subbundle is parallel with respect to
a metric linear connection with totally skew-symmetric torsion and the torsion
3-form is of type (1,2)+(2,1) with respect to all almost complex structures in $Q$. 
The target space of two-dimensional (4,0) supersymmetric sigma models with torsion 
coupled to (4,0) supergravity is a QKT manifold \cite{HOP}. 
If the torsion of a QKT manifold is a closed 3-form then it is 
called strong QKT manifold. The properties of HKT and
QKT geometries strongly resemble those of HK and QK ones, respectively. In
particular, HKT \cite{HP3} and QKT \cite{HOP} manifolds admit twistor constructions
with twistor spaces which have similar properties to those of HK \cite{HKLR} and QK
\cite{S1,S2,S3}. 

The main object of interest in this article is the differential geometric
properties of QKT manifolds. We find necessary and sufficient conditions to the
existence of a QKT connection in terms of the K\"ahler 2-forms and  show that the
QKT-connection is unique if dimension is at least 8 
(see Theorem~\ref{th1} below). We prove that the QKT manifolds
are invariant under  conformal transformations of the metric. This allows us to present
 a lot of 
(compact) examples of QKT manifolds. In particular, we show that the compact 
quaternionic Hopf manifolds studied in \cite{OP1}, which do not admit a QK structure, 
are QKT manifolds. In the compact case we show the existence of Gauduchon metric i.e.
the unique  conformally 
equivalent QKT structure with co-closed torsion 1-form.

It is shown in \cite{HOP} that the twistor space of a QKT manifold is always 
 complex manifold provided the dimension is at least 8.  
It admits complex contact (resp. K\"ahler) structure if the torsion 4-form is of
  type (2,2) and some additional nondegeneratity (positivity) conditions are  
 fulfilled \cite{HOP}. 
   Most of the known examples of QKT manifolds are homogeneous constructed 
   in\cite{homo}.  However, there are no 
homogeneous proper QKT manifolds (i.e. QKT which is not QK or HKT) 
with torsion 4-form of type (2,2) 
in dimensions greater than four by the result of \cite{homo}. 
We generalise this result showing that there are no proper QKT manifolds with 
torsion 4-form of type (2,2) provided 
that the torsion is parallel and dimension is at least 8. 

In dimension 4 a lot of examples of QKT manifolds are known \cite{HOP,homo}. 
In particular, 
examples of homogeneous QKT manifolds are constructed in  \cite{homo}.  
We notice that there are many (even strong) QKT structures in dimension 4,  
all depending on an arbitrary 1-form. 
We give a local description of 4-dimensional QKT manifolds with  
parallel torsion; namely such a QKT manifold is a Riemannian product of a real 
line and a 3-dimensional Riemannian manifold. 
We observe that homogeneous QKT manifolds are precisely naturally reductive 
homogeneous 
Riemannian manifolds, the objects which are well known. We present a complete 
local description 
(up to an isometry) of 4-dimensional homogeneous QKT which was known in the 
setting of naturally reductive homogeneous 4-manifold \cite{KV}. In the 
last section we consider 4-dimensional Einstein-like QKT manifold and 
find a closed relation with Einstein-Weyl geometry in dimension four. In particular,
we show that every 4-dimensional HKT manifold is of this type.

{\bf Acknowledgements.} The research was done during the author's visit 
at the Abdus Salam 
International Centre for Theoretical Physics, Trieste Italy. 
The author thanks the Abdus Salam ICTP for support and the 
excellent environment.
The author also thanks to  G. Papadopoulos for his interest, 
useful suggestions and remarks.  He is grateful to S.Marchiafava for pointing out
some incorrect statements 
and L.Ornea  for 
finding the time to read and comment on a first draft of the manuscript. 

\section{Characterisations of QKT connection}

Let $(M,g,(J_{\alpha})\in Q,\alpha=1,2,3)$ be a 4n-dimensional almost quaternionic
manifold with $Q$-hermitian Riemannian metric $g$ and an admissible basis
$(J_{\alpha})$. The K\"ahler form $F_{\alpha}$ of each $J_{\alpha}$ is defined by
$F_{\alpha}=g(.,J{\alpha}.)$. The corresponding Lee forms are given by
$\theta_{\alpha}=\delta F_{\alpha}\circ J_{\alpha}$. 

For an $r$-form $\psi$ we denote by $J_{\alpha}\psi$ the
$r$-form defined by\\
$J_{\alpha}\psi(X_1,...,X_r):=(-1)^r\psi(J_{\alpha}X_1,...,J_{\alpha}X_r),
\alpha=1,2,3$. Then $(d^c\psi)_{\alpha}=(-1)^rJ_{\alpha}d\psi$. We
shall use the notations 
$d_{\alpha}F_{\beta}:=(d^cF_{\beta})_{\alpha}$, i.e.
$d_{\alpha}F_{\beta}(X,Y,Z)=-dF_{\beta}(J_{\alpha}X,J_{\alpha}Y,J_{\alpha}Z), 
\alpha,\beta =1,2,3.$

We recall the decomposition of a skew-symmetric tensor $P\in \Lambda^2T^*M
\otimes TM$ with respect to a given almost complex structure $J_{\alpha}$. The
(1,1), (2,0) and (0,2) part of $P$ are defined by $P^{1,1}(J_{\alpha}X,J_{\alpha}Y)
= P^{1,1}(X,Y),  P^{2,0}(J_{\alpha}X,Y) = J_{\alpha}P^{2,0}(X,Y), 
 P^{0,2}(J_{\alpha}X,Y) = - J_{\alpha}P^{0,2}(X,Y)$, respectively.   

For each $\alpha = 1,2,3$, we denote by $dF^+_{\alpha}$ (resp. $dF^-_{\alpha}$) the
$(1,2)+(2,1)$-part (resp. $(3,0)+(0,3)$-part) of  $dF_{\alpha}$ with respect to the
almost complex structure $J_{\alpha}$. We consider the following 1-forms
$$
\theta_{\alpha,\beta} = -\frac{1}{2}\sum_{i=1}^{4n}dF^+_{\alpha}(X,e_i,J_{\beta}e_i), \quad \alpha ,\beta =1,2,3.
$$
Here and further $e_1,e_2,\dots,4n$ is an orthonormal basis of the tangential space. 

Note that $\theta_{\alpha,\alpha} = \theta_{\alpha}$.

The Nijenhuis tensor $N_{\alpha}$ of an almost complex structure $J_{\alpha}$ is
given by
\\
$N_{\alpha}(X,Y)=[J_{\alpha}X,J_{\alpha}Y]-[X,Y]-J_{\alpha}[J_{\alpha}X,Y] -
J_{\alpha}[X,J_{\alpha}Y].$

The celebrated Newlander-Nirenberg theorem \cite{NN} states that an
almost complex
structure is a complex structure if and only if its Nijenhuis tensor vanishes. 

Let $\nabla$ be a quaternionic connection i.e.
\begin{equation}\label{1}
\nabla J_{\alpha}=-\omega_{\beta}\otimes J_{\gamma} + \omega_{\gamma}\otimes 
J_{\beta},
\end{equation}
where the $\omega_{\alpha}, \alpha =1,2,3$ are 1-forms. 

Here and henceforth  
$(\alpha,\beta,\gamma)$
is a cyclic permutation of $(1,2,3)$.

Let $T(X,Y)=\nabla_XY-\nabla_YX-[X,Y]$ be the torsion tensor of type (1,2) of 
$\nabla$. We denote by the same letter the torsion tensor of type (0,3) given by
$T(X,Y,Z)=g(T(X,Y),Z)$. The Nijenhuis tensor is expressed in terms of $\nabla$ as
follows 
\begin{eqnarray}\label{2}
N_{\alpha}(X,Y) & = & 4T_{\alpha}^{0,2}(X,Y)\\
                & + &(\nabla_{J_{\alpha}X}J_{\alpha})(Y) -
(\nabla_{J_{\alpha}Y}J_{\alpha})(X) - (\nabla_Y J_{\alpha})(J_{\alpha}X) +
(\nabla_XJ_{\alpha})(J_{\alpha}Y),\nonumber
\end{eqnarray}
where the (0,2)-part $T_{\alpha}^{0,2}$ of the torsion with respect to
$J_{\alpha}$ is given by
\begin{equation}\label{tr1}
T_{\alpha}^{0,2}(X,Y) = \frac{1}{4}\left( T(X,Y) - T(J_{\alpha}X,J_{\alpha}Y) +
J_{\alpha}T(J_{\alpha}X,Y) + J_{\alpha}T(X,J_{\alpha}Y)\right).
\end{equation} 
We recall that if a 3-form $\psi$ is of type (1,2)+(2,1) with respect to an almost
complex structure $J$ then it satisfies the equality
\begin{equation}\label{3}
\psi(X,Y,Z)=\psi(JX,JY,Z)+\psi(JX,Y,JZ) + \psi(X,JY,JZ).
\end{equation}
{\bf Definition}. An almost quaternionic hermitian manifold $(M,g,(H_{\alpha})\in
Q)$ is {\it QKT manifold} if it admits a metric
quaternionic connection $\nabla$ with totally skew symmetric torsion which is
(1,2)+(2,1)-form with respect to each $J_{\alpha}, \alpha=1,2,3$.  If the torsion 
3-form is closed then the manifold is said to be {\it strong QKT manifold}.

It follows that the holonomy group of $\nabla$ is a subgroup of SP(n).SP(1). 

By means of (\ref{1}), (\ref{2}) and (\ref{3}), 
the Nijenhuis tensor $N_{\alpha}$ of $J_{\alpha},\alpha=1,2,3$, 
on a QKT manifold is given by 
\begin{equation}\label{6}
N_{\alpha}(X,Y)  = A_{\alpha}(Y)J_{\beta}X -  A_{\alpha}(X)J_{\beta}Y -
J_{\alpha} A_{\alpha}(Y)J_{\gamma}X + J_{\alpha} A_{\alpha}(X)J_{\gamma}Y,
\end{equation}
where
\begin{equation}\label{c1}
A_{\alpha}=\omega_{\beta} +J_{\alpha} \omega_{\gamma}.
\end{equation}
\noindent {\bf Remark 1.} The definition of QKT manifolds given above is equivalent
to that
given in
\cite{HOP} because the requirement the torsion to be (1,2)+(2,1)-form with respect
to each $J_{\alpha},\alpha=1,2,3$, is equivalent, by means of (\ref{6}), to the
fourth condition of (4) in \cite{HOP}. The torsion of $\nabla$ is (1,2)+(2,1)-form
with respect to any (local) almost complex structure $J\in Q$ \cite{HOP}. This
follows also from (\ref{6}) and the general formula (6) in \cite{AMP1} which expresses
$N_J$ in terms of  $N_{J_1}, N_{J_2}, N_{J_3}$. In fact, it is sufficient  that the
torsion
is a (1,2)+(2,1)-form with respect to the only two almost complex structures of
$(H)$ since the formula (3.4.4) in \cite{AM}) gives the necessary  
expression of $N_{J_3}$ by $N_{J_1}$ and $N_{J_2}$. Indeed, it is easy to see that
the formula (3.4.4) in \cite{AM} holds for the (0,2)-part $T_{\alpha}^{0,2}, \alpha
=1,2,3$, of the torsion. Hence, the vanishing of the (0,2)-part of the torsion
with respect to any two almost complex structures in $(H)$ implies the vanishing of
the (0,2)-part of $T$ with respect to the third one.

On a QKT manifold there are three naturally
associated 1-forms to the torsion  defined by
\begin{equation}\label{n1}
t_{\alpha}(X)=-\frac{1}{2}\sum_{i=1}^{4n}T(X,e_i,J_{\alpha}e_i), \quad \alpha=1,2,3.
\end{equation} 
We have
\begin{pro}\label{l1}
On a QKT manifold  $J_1t_1=J_2t_2=J_3t_3.$
\end{pro}
{\it Proof.} Applying (\ref{3}) with respect to $J_{\beta}$ we obtain 
\begin{eqnarray}
t_{\alpha}(X) &=& -\frac{1}{2}\sum_{i=1}^{4n} T(X,e_i,J_{\alpha}e_i) = 
 -\frac{1}{2}\sum_{i=1}^{4n}
T(X,J_{\beta}e_i,J_{\gamma}e_i) \nonumber\\
               &=& \frac{1}{2}\sum_{i=1}^{4n} T(J_{\beta}X,e_i,J_{\gamma}e_i)
	        -\frac{1}{2}\sum_{i=1}^{4n}
T(J_{\beta}X,J_{\beta}e_i,J_{\alpha}e_i) +\frac{1}{2} \sum_{i=1}^{4n}
T(X,e_i,J_{\alpha}e_i).\nonumber
\end{eqnarray}  
The last equality implies $t_{\alpha}=J_{\beta}t_{\gamma}$ which 
proves the assertion. \hfill {\bf Q.E.D.}

The 1-form $t=J_{\alpha}t_{\alpha}$ is independent of the chosen almost complex 
structure $J_{\alpha}$ by Proposition~\ref{l1}. We shall call it 
{\it the torsion 1-form} of a given QKT manifold.
\\[2mm] 
{\bf Remark 2.} Every QKT manifold is a quaternionic manifold. This is an
immediate consequence of (\ref{6}) and Proposition 2.3 in \cite{AMP1}.

However, the converse to the above property is not always true. In fact, we have
\begin{th}\label{th1}
Let $(M,g,(J_{\alpha}\in Q)$ be a 4n-dimensional ($n > 1$) quaternionic manifold
with $Q$-hermitian
metric $g$. Then $M$ admits a QKT structure if and only if the following conditions
hold
\begin{equation}\label{4}
(d_{\alpha}F_{\alpha})^+- (d_{\beta}F_{\beta})^+ 
=\frac{1}{2}\left( K_{\alpha}\wedge F_{\beta}- J_{\beta}K_{\beta}\wedge F_{\alpha}-
(K_{\beta}-J_{\alpha}K_{\alpha})\wedge F_{\gamma}\right),
\end{equation}
where $(d_{\alpha}F_{\alpha})^+$ denotes the (1,2)+(2,1) part of
$(d_{\alpha}F_{\alpha})$ with respect to the $J_{\alpha}, \alpha=1,2,3$. The 1-forms
$K_{\alpha}, \alpha=1,2,3$, are given by
\begin{equation}\label{c2}
K_{\alpha}= \frac{1}{1-n}\left(J_{\beta}\theta_{\alpha} + 
\theta_{\alpha,\gamma}\right).
\end{equation}
 The metric quaternionic
connection $\nabla$ with torsion
3-form of type (1,2)+(2,1) is unique and is determined  by
\begin{equation}\label{5} 
\nabla =\nabla^g +\frac{1}{2}\left((d_{\alpha}F_{\alpha})^+ -\frac{1}{2}
\left(J_{\alpha}K_{\alpha}\wedge F_{\gamma}+K_{\alpha}\wedge 
F_{\beta}\right)\right),
\end{equation}
where $\nabla^g$ is the Levi-Civita connection of $g$.
\end{th}
{\it Proof.} To prove the 'if' part, let $\nabla$ be a metric quaternionic connection
satisfying (\ref{1}) which torsion $T$ has the required properties. We follow the
scheme in \cite{Ga1}. Since $T$ is
skew-symmetric we have
\begin{equation}\label{5'}
\nabla =\nabla^g + \frac{1}{2}T.
\end{equation}
We obtain using (\ref{1}) and (\ref{5'}) that
\begin{eqnarray}\label{7}
\frac{1}{2}\left(T(X,J_{\alpha}Y,Z) + (T(X,Y,J_{\alpha}Z)\right) &=&
-g\left((\nabla_X^gJ_{\alpha})Y,Z\right)\\
&+&
\omega_{\beta}(X)F_{\gamma}(Y,Z) -   
\omega_{\gamma}(X)F_{\beta}(Y,Z).\nonumber
\end{eqnarray}
The tensor $\nabla^gJ_{\alpha}$ is decomposed by parts according to 
 $\nabla J_{\alpha} = 
(\nabla J_{\alpha})^{2,0} + (\nabla J_{\alpha})^{0,2}$, where \cite{Ga1}
\begin{equation}\label{10}
g\left((\nabla_X^gJ_{\alpha})^{2,0}Y,Z\right) =
\frac{1}{2}\left((dF_{\alpha})^+(X,J_{\alpha}Y,J_{\alpha}Z) -
(dF_{\alpha})^+(X,Y,Z)\right)
\end{equation}
\begin{equation}\label{9}
g\left((\nabla_X^gJ_{\alpha})^{0,2}Y,Z\right) = 
\frac{1}{2}\left(g(N_{\alpha}(X,Y),J_{\alpha}Z) - 
g(N_{\alpha}(X,Z),J_{\alpha}Y)-
g(N_{\alpha}(Y,Z),J_{\alpha}X)\right)
\end{equation}
Taking the (2,0) part in (\ref{7}) we obtain using (\ref{10}) that 
\begin{eqnarray}\label{i2}
T(X,J_{\alpha}Y,Z) + T(X,Y,J_{\alpha}Y) &=&
(dF_{\alpha}^+(X,J_{\alpha}Y,J_{\alpha}Z)-
(dF_{\alpha}^+(X,Y,Z)\\
&+&
C_{\alpha}(X)F_{\gamma}(Y,Z)+C_{\alpha}(J_{\alpha}X)F_{\beta}(Y,Z),\nonumber 
\end{eqnarray}
where 
\begin{equation}\label{c3}
C_{\alpha}=\omega_{\beta}-J_{\alpha}\omega_{\gamma}.
\end{equation}
The cyclic sum of (\ref{i2}) and  the fact that $T$ and
$(dF_{\alpha})^+$
are (1,2)+(2,1)-forms with respect to each $J_{\alpha}$, gives
\begin{equation}\label{i3}
T= (d_{\alpha}F_{\alpha})^+ -\frac{1}{2}
\left(J_{\alpha}C_{\alpha}\wedge F_{\gamma}+C_{\alpha}\wedge
F_{\beta}\right).
\end{equation}
Further, we take the contractions in (\ref{i3}) to get
\begin{eqnarray}\label{c4} 
J_{\alpha}t_{\alpha} = 
- \theta_{\alpha} - J_{\beta}C_{\alpha}, \nonumber\\
J_{\alpha}t_{\alpha} = 
- J_{\gamma}\theta_{\beta,\alpha} - n J_{\gamma}C_{\beta},\\
J_{\alpha}t_{\alpha} = 
J_{\beta} \theta_{\gamma,\alpha} - n J_{\alpha}C_{\gamma}\nonumber
\end{eqnarray}
Using Proposition~\ref{l1}, (\ref{c1}) and (\ref{c3}), we obtain 
consequently from  
(\ref{c4})  that
\begin{equation}\label{c5}
A_{\alpha}= J_{\alpha}C_{\beta} +J_{\gamma}C_{\gamma} = 
J_{\beta}\left(\theta_{\gamma} - \theta_{\beta}\right),
\end{equation}
\begin{equation}\label{c6}
(n-1)J_{\beta}C_{\alpha}= \theta_{\alpha} - 
J_{\beta}\theta_{\alpha,\gamma}. 
\end{equation}
Then (\ref{4}) and (\ref{c2}) follow from (\ref{i3}) and (\ref{c6}). 
                                         
For the converse, we define $\nabla$ by (\ref{5}). To complete the proof we have to
show that $\nabla$ is a quaternionic connection. We calculate
\begin{eqnarray}
g\left((\nabla_XJ_{\alpha})Y,Z\right) &=&
g\left((\nabla_X^gJ_{\alpha})Y,Z\right)+
\frac{1}{2}\left(T(X,J_{\alpha}Y,Z) +
T(X,Y,J_{\alpha}Z)\right) \nonumber\\
                                      &=&
\omega_{\beta}(X)F_{\gamma}(Y,Z) -   \omega_{\gamma}(X)F_{\beta}(Y,Z),\nonumber
\end{eqnarray}
where we used (\ref{10}), (\ref{9}), (\ref{c5}), (\ref{c2}), (\ref{c1}), (\ref{c3})
and the compatibility condition (\ref{4}) to get the last equality. 
The uniqueness of $\nabla $ follows from (\ref{5}) as well as from Theorem 10.3 in 
\cite{Ob} which states that any quaternionic connection is entirely determined 
by its torsion (see also \cite{Gau}).
\hfill {\bf Q.E.D.}

In the case of HKT manifold, $K_{\alpha}=dF_{\alpha}^-=0$ and Theorem~\ref{th1} 
is a consequence of the general results in \cite{Ga1} (see also \cite{GP}) 
which imply  
that on a hermitian manifold there exists a unique linear connection with totally
skew-symmetric torsion preserving the metric and the complex structure, the Bismut 
connection. This connection was used by Bismut \cite{BI} to prove a local index 
theorem for the Dolbeault operator on non-K\"ahler manifold.  The geometry of 
this connection is referred to KT-geometry by physicists. Obstructions to the 
existence of (non-trivial) Dolbeault cohomology groups on a compact KT-manifold 
are presented in \cite{AI}.

We note that (\ref{c5}) and (\ref{c6}) are also valid in the case $n=1$.
 
We get, as a consequence of the proof of Theorem~\ref{th1},  the following
integrability criterion which is discovered in dimension 4 in \cite{GT}. 
\begin{pro}\label{cc1}
The Nijenhuis tensors  of a QKT manifold depend only on the difference between 
the Lie forms. In particular, 
the almost complex structures $J_{\alpha}$ on a QKT manifold $(M,(J_{\alpha})\in
Q,g,\nabla )$ are integrable if and only if
$$\theta_{\alpha}=\theta_{\beta}= \theta_{\gamma}$$
\end{pro}
{\it Proof.} The Nijenhuis tensors  are given by (\ref{6}) and (\ref{c5}). 
\hfill {\bf Q.E.D.}

\begin{co}\label{tt1}
On a 4n-dimensional QKT manifold the following formulas hold
$$J_{\beta}\theta_{\alpha,\gamma}=
- J_{\gamma}\theta_{\alpha,\beta},$$
\begin{equation}\label{per1}
(n^2+n)\theta_{\alpha} -n\theta_{\beta} - n^2\theta_{\gamma} +
J_{\gamma}\theta_{\beta,\alpha} + nJ_{\alpha}\theta_{\gamma,\beta} -
(n+1)J_{\beta}\theta_{\alpha,\gamma}=0.
\end{equation}
If $n=1$ then  \hspace{4cm} $\theta_{\alpha} = J_{\beta}
\theta_{\alpha,\gamma} = 
-J_{\gamma}\theta_{\alpha,\beta}$.
\end{co}
{\it Proof.} The first formula follows directly from the system 
(\ref{c4}). Solving the system (\ref{c4}) with respect to $C_{\alpha}$ we obtain
\begin{equation}\label{per2}
(n^3-1)J_{\beta}C_{\alpha} = (\theta_{\alpha}-J_{\gamma}\theta_{\beta,\alpha}) +
n(\theta_{\beta}-J_{\alpha}\theta_{\gamma,\beta}) +
n^2(\theta_{\gamma}-J_{\beta}\theta_{\alpha,\gamma}).
\end{equation}
Then (\ref{per1}) is a consequence of (\ref{per2}) and (\ref{c6}). 
 The last assertion follows from (\ref{c6}) . 
\hfill {\bf Q.E.D.}
\begin{co}\label{cc2}
On a 4n-dimensional ($n>1$) QKT manifold the $sp(1)$-connection 1-forms are given by
\begin{equation}\label{c7}
\omega_{\beta} = \frac{1}{2}J_{\beta}\left(\theta_{\gamma} - \theta_{\beta}
+ \frac{1}{1-n}\theta_{\alpha}\right)  
+ \frac{1}{2(1-n)}\theta_{\alpha,\gamma}. 
\end{equation}
\end{co}
{\it Proof.} The proof follows in a straightforward way from (\ref{c5}), 
(\ref{c6}), (\ref{c1}) and (\ref{c3}). \hfill {\bf Q.E.D.}

Theorem~\ref{th1} and the above formulas lead to the following criterion
\begin{pro}\label{cc3}
Let $(M,g,(H))$ be a 4n-dimensional ($n>1$) QKT manifold. The following conditions are equivalent:

i) $(M,g,(H))$ is a HKT manifold;

ii) $d_{\alpha}F_{\alpha}^+=d_{\beta}F_{\beta}^+= d_{\gamma}F_{\gamma}^+$;

iii) 
$
\theta_{\alpha} = J_{\beta}\theta_{\gamma,\alpha}.
$
\end{pro}
{\it Proof.} 
If $(M,g,(H))$ is a HKT manifold, the connection 1-forms 
$\omega_{\alpha}=0, \alpha=1,2,3$. 
Then ii) and iii) follow from (\ref{c3}), (\ref{c6}), 
(\ref{c2}) and (\ref{4}).

If iii) holds, then (\ref{c6}) and (\ref{c5}) yield $C_{\alpha}=A_{\alpha}=0, 
\alpha=1,2,3$, since $n>1$. Consequently, $2\omega_{\alpha}= 
J_{\beta}C_{\beta}- J_{\beta}A_{\beta}=0$ by (\ref{c3}) and 
(\ref{c1}). Thus the equivalence of i) and iii) is proved.

Let ii) holds. Then we compute that $\theta_{\alpha}=
J_{\gamma}\theta_{\beta,\alpha}$. Since $n>1$, the equality (\ref{per2}) 
leads to $C_{\alpha}=0,
\alpha=1,2,3$, which forces $\omega_{\alpha}=0, \alpha=1,2,3$  as above. 
This completes the proof. \hfill {\bf Q.E.D.}

The next theorem shows that QKT manifolds are stable under a conformal transformations.
\begin{th}\label{th3}
Let $(M,g,(J_{\alpha}),\nabla)$ be a 4n-dimensional QKT manifold. Then every Riemannian 
metric $\bar{g}$ in the conformal class [g] admits a QKT connection. If $\bar{g}=fg$ for 
a positive function $f$ then the QKT connection $\bar{\nabla}$ corresponding to $\bar{g}$
 is given by
\begin{eqnarray}\label{z1}
\bar{g}(\bar{\nabla}_XY,Z)=fg(\nabla_XY,Z) &+&
\frac{1}{2}\left(df(X)g(Y,Z) + df(Y)g(X,Z) - df(Z)g(X,Y)\right)\\
&+&
\frac{1}{2}\left(J_{\alpha}df\wedge F_{\alpha} +J_{\beta}df\wedge F_{\beta}+
J_{\gamma}df\wedge F_{\gamma}\right)(X,Y,Z).\nonumber
\end{eqnarray}

If $M$ is compact then there exists a unique (up to homotety) metric $g_G \in [g]$
with co-closed torsion 1-form.
\end{th}
{\it Proof.} First we assume $n>1$. We shall apply Theorem~\ref{th1} to the quaternionic
Hermitian manifold $(M,\bar{g}=fg,(J_{\alpha})\in Q)$. We denote the objects corresponding 
to the metric $\bar{g}$ by a line above the symbol e.g. $\bar{F_{\alpha}}$ denotes the K\"ahler 
form of $J_{\alpha}$ with respect to $\bar{g}$. An easy calculation
gives the following
sequence of formulas
\begin{equation}\label{z2}
d_{\alpha}\bar{F}_{\alpha}^+ = J_{\alpha}df\wedge F_{\alpha} + fd_{\alpha}F_{\alpha}^+;
\quad 
\bar{\theta}_{\alpha} = \theta_{\alpha} +(2n-1)d\ln f; \quad 
\bar{\theta}_{\alpha,\gamma} = \theta_{\alpha,\gamma} - J_{\beta}d\ln f.
\end{equation}
We substitute (\ref{z2}) into (\ref{c2}), (\ref{c5}) and (\ref{c7}) to get
\begin{equation}\label{z3}
\bar{K}_{\alpha} = K_{\alpha} -2J_{\beta}d\ln f, \quad \bar{A}=A, \quad
\bar{\omega}_{\alpha} = \omega_{\alpha} -J_{\beta}d\ln f.
\end{equation}
Using (\ref{z2}) and (\ref{z3}) we verify that the conditions (\ref{4}) 
with respect to the metric $\bar{g}$ are fulfilled. Theorem~\ref{th1} implies that 
there exists a QKT connection $\bar{\nabla}$ with respect to  $(\bar{g},Q)$. 
Using the well known relation between the Levi-Civita connections of conformally 
equivalent metrics, (\ref{z2}) and (\ref{z3}), we obtain  
(\ref{z1}) from (\ref{5}).

If $n=1$ we define the new QKT connection with respect to $(\bar{g},Q)$  by (\ref{z1}).

Using (\ref{z1}), we find that the torsion tensors $T$ and $\bar{T}$ of $\nabla$ 
and $\bar{\nabla}$ 
are related by
\begin{equation}\label{z4}
\bar{T}= fT + J_{\alpha}df \wedge F_{\alpha} +J_{\beta}df \wedge F_{\beta}+
J_{\gamma}df \wedge F_{\gamma}.
\end{equation}
Consequently, we obtain from (\ref{z4}) for the torsion 1-forms $t$ and $\bar{t}$ 
that 
\begin{equation}\label{z5}
\bar{t}= t - (2n+1)d\ln f.
\end{equation}
If $M$ is compact, we may apply to (\ref{z5}) the theorem of Gauduchon for the existence of a 
Gauduchon metric on a compact Weyl manifold \cite{G1,G2} to obtain the desired metric 
$g_G$. \hfill {\bf Q.E.D.} 

We shall call the unique metric with co-closed torsion 1-form on a compact QKT 
manifold the {\it Gauduchon metric}.
\begin{co}
On a compact QKT manifold with closed (non exact) torsion 1-form the Gauduchon 
metric $g_G$ cannot have positive definite Riemannian Ricci tensor. In particular, 
if it is an Einstein manifold then it is of non-positive scalar curvature. 

Further, if the Gauduchon metric is Ricci flat then the corresponding  torsion 
1-form $t_G$ is parallel with respect to the Levi-Civita connection of $g_G$.
\end{co}
{\it Proof.} The two form $dt$ is invariant under conformal transformations by 
(\ref{z5}). Then the Gauduchon metric has harmonic torsion 1-form i.e 
$dt= \delta t=0$. 
The claim follows from the  Weitzenb\"oeck formula (see e.g. \cite{Bes})  
$
\int_{M}\{|dt|^2 +|\delta t|^2\} \,dV = \int_{M} \{|\nabla^gt|^2+
Ric^g(t^{\#},t^{\#})\} \,dV=0,$ 
where $t^{\#}$ is 
the dual vector field of $t$, $|.|$ is the usual tensor norm and $dV$ is the volume 
form. \hfill {\bf Q.E.D.}

Theorem~\ref{th3} allows us to supply a large class of (compact) QKT manifold. 
Namely, any conformal metric of a QK, HK or HKT manifold will give a QKT manifold. 
This leads to the notion of {\it locally conformally QK (resp. locally conformally HK, resp.
locally conformally HKT) 
manifolds} 
(briefly l.c.QK (resp. l.c.HK, resp. l.c.HKT) manifolds) in the context of QKT geometry.

The l.c.QK and l.c.HK manifolds have already appeared in the context of Hermitian-Einstein-Weyl
structures \cite{PPS} and of 3-Sasakian structures \cite{BGM}. These two classes of  
quaternionic 
manifolds are studied in detail (mostly in the compact case) in \cite{OP1,OP2}. 

We recall that a quaternionic Hermitian manifold $(M,g,Q)$ is said to be  l.c.QK (resp.
l.c.HK, resp. l.c.HKT) manifold if each point $p\in M$ has a neighbourhood $U_p$ such that 
$g \Big|_{U_p}$ is conformally equivalent to a QK (resp.HK, resp.HKT) metric. 
There are compact l.c.QK manifold which do not admit any QK structure \cite{OP1}.
Typical examples of compact l.c. QK manifolds without any QK structure are the quaternionic 
Hopf spaces $H=({\cal H}^n-\{0\})/\Gamma$, where $\Gamma$ is an appropriate discrete 
group acting diagonally on the quaternionic coordinates in ${\cal H}^n$ (see \cite{OP1}).

We recall that on a l.c.QK manifold the 4-form $\Omega =\sum_{\alpha=1}^{3} 
F_{\alpha}\wedge F_{\alpha}$ satisfies $d\Omega=\omega\wedge \Omega, d\omega=0$, 
where $\omega $ is locally defined by $\omega=2d\ln f$. On a l.c.QK manifold viewed 
as a QKT manifold by Theorem~\ref{th3} the torsion 1-form is equal to 
$t=(2n+1)\omega$ by (\ref{z5}). The QK manifolds are Einstein provided the 
dimension is at least 8 \cite{A1,Ber}. Then, the Gauduchon Theorem \cite{G2} applied to l.c.QK 
manifold in \cite{OP1} can be stated in our context as follows
\begin{co}
Let $(M,g)$ be a compact 4n-dimensional ($n>1$) QKT manifold which is l.c.QK and assume that 
no metric in the conformal class [g] of g is QK. Then the torsion 1-form of the 
Gauduchon metric $g_G$ is parallel with respect to the Levi-Civita connection of $g_G$.
\end{co}
Theorem~\ref{th3},  Theorem~\ref{th1} together with  Proposition~\ref{cc1}
and Proposition~\ref{cc3} imply  the following
\begin{co}\label{u1}
Every l.c.QK manifold admits a QKT structure. 

Further, if  $(M,g,(J_{\alpha}),\nabla)$ is a 4n-dimensional $n>1$ QKT manifold then: 

i) $(M,g,(J_{\alpha}),\nabla)$ is a l.c.QK manifold if and only if 
\begin{equation}\label{u3}
T = \frac{1}{2n+1}\left(t_{\alpha}\wedge F_{\alpha} +
t_{\beta}\wedge 
F_{\beta}+t_{\gamma}\wedge F_{\gamma}\right), \quad dt=0;
\end{equation}
\indent ii) $(M,g,(J_{\alpha}),\nabla)$ is a l.c.HKT manifold if and only if 
the 1-form $\theta_{\alpha}-J_{\beta}\theta_{\alpha,\gamma}$ is closed i.e. 
$$
d(\theta_{\alpha}-J_{\beta}\theta_{\alpha,\gamma})=0;
$$
\indent iii)  $(M,g,(J_{\alpha}),\nabla)$ is a l.c.HK manifold if an only if (\ref{u3})
holds and 
$$
\theta_{\alpha}-J_{\beta}\theta_{\alpha,\gamma}=\frac{2(1-n)}{2n+1}t.
$$
\end{co}

\section{Curvature of a QKT space}

Let $R=[\nabla,\nabla]-\nabla_[,]$ be the curvature tensor of type (1,3) of
$\nabla$. We denote the curvature tensor of type (0,4) $R(X,Y,Z,V)=g(R(X,Y)Z,V)$ by
the same letter. There are three Ricci forms given by
$$
\rho_{\alpha}(X,Y)=\frac{1}{2} \sum_{i=1}^{4n}R(X,Y,e_i,J_{\alpha}e_i), \quad
\alpha=1,2,3.
$$
\begin{pro}\label{p1}
The curvature of a QKT manifold $(M,g,(J_{\alpha}),\nabla)$ 
satisfies the following
relations
\begin{equation}\label{11}
R(X,Y)J_{\alpha} = \frac{1}{n}\left(\rho_{\gamma}(X,Y)J_{\beta} 
-\rho_{\beta}(X,Y)J_{\gamma}\right),
\end{equation}
\begin{equation}\label{12}
\rho_{\alpha} =d\omega_{\alpha} + \omega_{\beta}\wedge\omega_{\gamma}.
\end{equation}
\end{pro}
{\it Proof.} We follow the classical scheme (see e.g. \cite{AM,Ish,Bes}). Using (\ref{1})
we obtain
$$
R(X,Y)J_{\alpha}=-(d\omega_{\beta} +
\omega_{\gamma}\wedge\omega_{\alpha})(X,Y)J_{\gamma} + 
(d\omega_{\gamma} +
\omega_{\alpha}\wedge\omega_{\beta})(X,Y)J_{\beta}.
$$
Taking the trace in the last equality, we get
\begin{eqnarray}
\rho_{\alpha}(X,Y) &=& 
\frac{1}{2} \sum_{i=1}^{4n}R(X,Y,e_i,J_{\alpha}e_i)=
\frac{1}{2} \sum_{i=1}^{4n}R(X,Y,J_{\beta}e_i,J_{\gamma}e_i) \nonumber\\
                   &=&
-\frac{1}{2} \sum_{i=1}^{4n}R(X,Y,e_i,J_{\alpha}e_i) 
+2n(d\omega_{\alpha} +
\omega_{\beta}\wedge\omega_{\gamma})(X,Y)J_{\beta}.\nonumber
\end{eqnarray}
\hfill {\bf Q.E.D.}

Using Proposition~\ref{p1} we find a simple necessary and sufficient condition a QKT
manifold to be a HKT one, i.e. the holonomy group of $\nabla$ to be a subgroup of
Sp(n).
\begin{pro}\label{p2}
A 4n-dimensional $(n>1)$ QKT manifold is a HKT manifold if and only if all the three Ricci forms vanish,
i.e $\rho_1=\rho_2=\rho_3=0$.
\end{pro}
{\it Proof.} If a QKT manifold is a HKT manifold then the holonomy group of $\nabla$ is
contained in Sp(n). This implies $\rho_{\alpha}=0, \quad \alpha=1,2,3$.

For the converse, let the three Ricci forms vanish. The equations (\ref{12}) mean
that the curvature of the Sp(1) connection on $Q$ vanish. Then there exists a
basis $(I_{\alpha}, \alpha =1,2,3)$ of almost complex structures on $Q$ and each 
$I_{\alpha}$ is $\nabla$-parallel i.e. the corresponding connection 1-forms 
$\omega_{I_{\alpha}}=0, \alpha =1,2,3$. Then each $I_{\alpha}$ is a complex 
structure, by
(\ref{6}) and (\ref{c1}). This implies that the QKT manifold is a HKT manifold. 
\hfill {\bf Q.E.D.}

We denote by $Ric,Ric^g$ the Ricci tensors of the QKT connection and of the
Levi-Civita
connection, respectively. In fact $Ric(X,Y)=\sum_{i=1}^{4n}R(e_i,X,Y,e_i)$.

Our main technical result is the following
\begin{pro}\label{tir}
Let $(M,g,(J_{\alpha}),\nabla)$ be a 4n-dimensional QKT manifold. The following 
formulas hold 
\begin{eqnarray}\label{ti20}
& &n\rho_{\alpha}(X,J_{\alpha}Y)+ \rho_{\beta}(X,J_{\beta}Y) +
\rho_{\gamma}(X,J_{\gamma}Y)= \\ 
& &
- nRic(XY)
+\frac{n}{4}(dT)_{\alpha}(X,J_{\alpha}Y)+ 
\frac{n}{2}(\nabla T)_{\alpha}(X,J_{\alpha}Y);\nonumber
\end{eqnarray}
\begin{eqnarray}\label{22}
& &(n-1)\rho_{\alpha}(X,J_{\alpha}Y)=-\frac{n(n-1)}{n+2}Ric(X,Y) \\
& &
+ \frac{n}{4(n+2)}\left\{(n+1)(dT)_{\alpha}(X,J_{\alpha}Y)-
(dT)_{\beta}(X,J_{\beta}Y)
-(dT)_{\gamma}(X,J_{\gamma}Y)\right\}\nonumber \\
& &
+ \frac{n}{2(n+2)}\left\{(n+1)(\nabla T)_{\alpha}(X,J_{\alpha}Y)-
(\nabla T)_{\beta}(X,J_{\beta}Y)
-(\nabla T)_{\gamma}(X,J_{\gamma}Y)\right\},
\end{eqnarray}
where

  $(dT)_{\alpha}(X,Y)=\sum_{i=1}^{4n}dT(X,Y,e_i,J_{\alpha}e_i), \quad 
 (\nabla T)_{\alpha}(X,Y) =\sum_{i=1}^{4n}(\nabla_XT)(Y,e_i,J_{\alpha}e_i)$.
\end{pro}
{\it Proof.} Since the torsion is a 3-form, we have
\begin{equation}\label{sof}
(\nabla^g_XT)(Y,Z,U) = (\nabla_XT)(Y,Z,U) + \frac{1}{2}
{\sigma \atop XYZ}
\left\{g(T(X,Y),T(Z,U)\right\}, 
\end{equation}
where ${\sigma \atop XYZ}$ denote the cyclic sum of $X,Y,Z$.
 
The exterior derivative $dT$ is 
given by
\begin{eqnarray}\label{13}
dT(X,Y,Z,U)&=&
{\sigma \atop XYZ}\left\{(\nabla_XT)(Y,Z,U) + g(T(X,Y),T(Z,U)\right\}\\
           &-& (\nabla_UT)(X,Y,Z) + {\sigma \atop XYZ}
\left\{g(T(X,Y),T(Z,U)\right\}. \nonumber
\end{eqnarray}

The first Bianchi identity for $\nabla$ states
\begin{equation}\label{14}
{\sigma \atop XYZ}R(X,Y,Z,U)= {\sigma \atop XYZ}\left\{(\nabla_XT)(Y,Z,U) +
g(T(X,Y),T(Z,U)\right\}.
\end{equation}
We denote by $B$ the Bianchi projector i.e. $B(X,Y,Z,U)={\sigma \atop
XYZ}R(X,Y,Z,U)$.

The curvature $R^g$ of the Levi-Civita connection is connected by $R$ in the
following way
\begin{eqnarray}\label{15}
R^g(X,Y,Z,U) &=& R(X,Y,Z,U) - \frac{1}{2} (\nabla_XT)Y,Z,U)
+\frac{1}{2} (\nabla_YT)X,Z,U)\nonumber \\
             &-&
\frac{1}{2}g(T(X,Y),T(Z,U)) - \frac{1}{4}g(T(Y,Z),T(X,U)) -
\frac{1}{4}g(T(Z,X),T(Y,U)).
\end{eqnarray}
Define $D$ by $D(X,Y,Z,U)=R(X,Y,Z,U)-R(Z,U,X,Y)$, we obtain from (\ref{15})
\begin{eqnarray}\label{tir1}
& &D(X,Y,Z,U) =\\
& &
\frac{1}{2}(\nabla_XT)(Y,Z,U)-\frac{1}{2}(\nabla_YT)(X,Z,U) -
\frac{1}{2}(\nabla_ZT)(U,X,Y)+ \frac{1}{2}(\nabla_UT)(Z,X,Y),\nonumber
\end{eqnarray}
since $D^g$ of $R^g$ is zero.

Using (\ref{11}) and (\ref{14}) we find the following relation between the Ricci
tensor and the Ricci forms
\begin{eqnarray}\label{16}
\rho_{\alpha}(X,Y) &=&
-\frac{1}{2} \sum_{i=1}^{4n}\left(R(Y,e_i,X,J_{\alpha}e_i)+
R(e_i,X,Y,J_{\alpha}e_i)\right) +\frac{1}{2}
\sum_{i=1}^{4n}B(X,Y,e_i,J_{\alpha}e_i)\\
                   &=& 
-\frac{1}{2}Ric(Y,J_{\alpha}X) +\frac{1}{2}Ric(X,J_{\alpha}Y) 
+\frac{1}{2}\sum_{i=1}^{4n}B(X,Y,e_i,J_{\alpha}e_i) \nonumber\\ 
                   &+&
\frac{1}{2n}\left\{\rho_{\beta}(J_{\gamma}Y,X)
-\rho_{\beta}(J_{\gamma}X,Y)
+\rho_{\gamma}(J_{\beta}X,Y)-\rho_{\gamma}(J_{\beta}Y,X)\right\}. \nonumber
\end{eqnarray}

On the other hand,  using  (\ref{11}), we calculate
\begin{eqnarray}\label{19}
\sum_{i=1}^4D(X,e_i,J_{\alpha}e_i,Y) &=&\sum_{i=1}^{4n}\{R(X,e_i,J_{\alpha}e_i,Y)+
R(Y,e_i,J_{\alpha}e_iX)\}\\ 
  &=&-Ric(Y,J_{\alpha}X)-Ric(X,J_{\alpha}Y)\nonumber \\
  &+&
\frac{1}{n}\left\{\rho_{\beta}(X,J_{\gamma}Y)
+\rho_{\beta}(Y,J_{\gamma}X)
-\rho_{\gamma}(Y,J_{\beta}X)-\rho_{\gamma}(X,J_{\beta}Y)\right\}. \nonumber
\end{eqnarray}
Combining (\ref{16}) and (\ref{19}), we derive
\begin{equation}\label{20}
n\rho_{\alpha}(X,J_{\alpha}Y)+ \rho_{\beta}(X,J_{\beta}Y) +
\rho_{\gamma}(X,J_{\gamma}Y)=
\end{equation}
$$
 - nRic(XY)+\frac{n}{2}B_{\alpha}(X,J_{\alpha}Y)+ 
\frac{n}{2}D_{\alpha}(X,J_{\alpha}Y), 
$$
where the tensors $B_{\alpha}$ and $D_{\alpha}$  are  defined by
$B_{\alpha}(X,Y)=\sum_{i=1}^{4n}B(X,Y,e_i,J_{\alpha}e_i)$ and \\
$
D_{\alpha}(X,Y)=\sum_{i=1}^{4n}D(X,e_i,J_{\alpha}e_i,Y)
$. Taking into account (\ref{tir1}), we get the expression
\begin{equation}\label{tir2}
D_{\alpha}(X,Y)=
\frac{1}{2}\sum_{i=1}^{4n}(\nabla_XT)(Y,e_i,J_{\alpha}e_i) +
\frac{1}{2}\sum_{i=1}^{4n}(\nabla_YT)(X,e_i,J_{\alpha}e_i)\quad \alpha=1,2,3.
\end{equation}
To calculate $B_{\alpha}+D_{\alpha}$ we use (\ref{13}) twice and (\ref{tir2}). 
After some calculations, we derive
\begin{equation}\label{tir4}
B_{\alpha}(X,Y)+ D_{\alpha}(X,Y) =\frac{1}{2}  \sum_{i=1}^{4n}
dT(X,Y,e_i,J_{\alpha}e_i) + \sum_{i=1}^{4n}(\nabla_XT)(Y,e_i,J_{\alpha}e_i), 
\quad \alpha=1,2,3.
\end{equation}
We substitute (\ref{tir4}) into (\ref{20}). Solving  the obtained system, we 
obtain
\begin{eqnarray}\label{21}
& &(n-1)\left\{\rho_{\alpha}(X,J_{\alpha}Y)- \rho_{\beta}(X,J_{\beta}Y)\right\} =\\
& &
\frac{n}{2}\left\{(dT)_{\alpha}(X,J_{\alpha}Y) - (dT)_{\beta}(X,J_{\beta}Y)\right\} +
\frac{n}{2}\left\{(\nabla T)_{\alpha}(X,J_{\alpha}Y) - 
(\nabla T)_{\beta}(X,J_{\beta}Y)\right\}.\nonumber
\end{eqnarray}
Finally, (\ref{20}) and (\ref{21}) imply (\ref{22}). \hfill {\bf Q.E.D.}
\\[2mm]
{\bf Remark 3.} The Ricci tensor of a QKT connection is not symmetric in general.
From (\ref{14}), (\ref{sof}) and the fact that $T$ is a 3-form we get the formula 
$
Ric(X,Y)-Ric(Y,X) = \sum_{i=1}^{4n}(\nabla^g_{e_i}T)(e_i,X,Y)=-\delta T(X,Y).
$
Hence, the Ricci tensor of a linear connection with totally skew-symmetric 
torsion is symmetric if and only if the torsion 3-form is co-closed.

\section{QKT manifolds with parallel torsion and homogeneous QKT structures}

Let $(G/K,g)$ be a reductive (locally) homogeneous Riemannian manifold. The canonical connection 
$\nabla$ is characterised by the properties $\nabla g=\nabla T=\nabla R=0$ 
\cite{KN},p.193.
A  homogeneous quaternionic Hermitian manifold (resp. homogeneous hyper Hermitian) 
manifold $(G/K,g,Q)$ is a homogeneous Riemannian manifold with an invariant quaternionic
Hermitian subbundle $Q$ (resp.  three invariant anti commuting complex structures ). 
This means 
that the bundle $Q$ (resp. each of the three complex structures) is parallel with respect 
to the canonical connection $\nabla$. The 
torsion of $\nabla$ is totally skew-symmetric if and only if the homogeneous 
Riemannian manifold is naturally reductive \cite{KN} (see also \cite{TV,homo}. 
Homogeneous QKT (resp. HKT) manifolds are homogeneous quaternionic 
Hermitian (resp. homogeneous hyper Hermitian) manifold which are naturally reductive.  
Examples of homogeneous HKT and QKT manifolds are presented in \cite{homo}. The 
homogeneous QKT manifolds in \cite{homo} are constructed from homogeneous HKT 
manifolds.

In this section we generalise the result of \cite{homo} which states that there are 
no homogeneous QKT manifold with torsion 4-form $dT$ of type (2,2) 
in dimensions greater than four. First, we prove the following technical result
\begin{pro}\label{teh}
Let $(M,g,(J_{\alpha}),\nabla)$ be a 4n-dimensional $(n>1)$ QKT manifold with 
4-form $dT$ 
of type (2,2) with respect to each $J_{\alpha}, \alpha =1,2,3$. Suppose that the 
torsion 
is parallel with respect to the QKT-connection. Then the Ricci forms $\rho_{\alpha}$ 
 are given by
\begin{equation}\label{27}
\rho_{\alpha}(X,J_{\alpha}Z) = \lambda g(X,Y), \quad \alpha =1,2,3,
\end{equation}
where $\lambda$ is a smooth function on $M$.
\end{pro}
{\it Proof.} 
Let the torsion be parallel i.e. $\nabla T=0$. Remark 3 shows that the Ricci tensor
is symmetric. The equalities (\ref{13}) and
(\ref{14}) imply
\begin{equation}\label{17}
B(X,Y,Z,U)=  {\sigma \atop XYZ}
\left\{g(T(X,Y),T(Z,U)\right\}=\frac{1}{2}dT(X,Y,Z,U).
\end{equation}
We get $D=0$ from (\ref{tir1}).

Suppose now that the 4-form $dT$ is of type (2,2) with respect to each $J_{\alpha},
\alpha
=1,2,3.$. Then it satisfies the  equalities
\begin{equation}\label{23}
dT(X,Y,Z,U)=dT(J_{\alpha}X,J_{\alpha}Y,Z,U) +
dT(J_{\alpha}X,Y,J_{\alpha}Z,U) + dT(X,J_{\alpha}Y,J_{\alpha}Z,U).
\end{equation}
The similar arguments as we used in the proof of Proposition~\ref{l1} but 
 applying 
(\ref{23}) instead of (\ref{3}), yield
\begin{lem}\label{l2}
On a QKT manifold with 4-form $dT$ of type (2,2) with respect to each $J_{\alpha},
\alpha=1,2,3$, the following equalities hold:
\begin{equation}\label{24}
(dT)_1(X,J_1Y) = (dT)_2(X,J_2Y) = (dT)_3(X,J_3Y), 
\end{equation}
\begin{equation}\label{24'}
(dT)_{\alpha}(X,J_{\alpha}Y)= -(dT)_{\alpha}(J_{\alpha}X,Y), \quad \alpha=1,2,3.
\end{equation}
\end{lem}
We substitute (\ref{24}), (\ref{17}) and $D=0$ into (\ref{21}) and (\ref{22}) 
to get
\begin{equation}\label{a1}
\rho_1(X,J_1Y) = \rho_2(X,J_2Y) = \rho_3(X,J_3Y),
\end{equation}
\begin{equation}\label{25}
\rho_{\alpha}(X,J_{\alpha}Y)=-\frac{n}{n+2}Ric(X,Y)
+ \frac{n}{4(n+2)}(dT)_{\alpha}(X,J_{\alpha}Y), \quad \alpha=1,2,3.
\end{equation}
The equality (\ref{24'}) shows that the 2-form $dT_{\alpha}$ is a (1,1)-form 
with respect
to $J_{\alpha}$. Hence, the $dT_{\alpha}$ is (1,1)-form with respect to each
$J_{\alpha},\alpha=1,2,3$, because of (\ref{24}). Since the Ricci tensor $Ric$ 
is symmetric, (\ref{25}) shows 
 that the Ricci tensor  $Ric$ is of hybrid type with respect to each
$J_{\alpha}$ i.e.
$Ric(J_{\alpha}X,J_{\alpha}Y)=Ric(X,Y), \alpha =1,2,3$ and  the Ricci forms
$\rho_{\alpha}, \alpha =1,2,3$ are  (1,1)-forms with respect to all $J_{\alpha},
\alpha=1,2,3$.
Taking into account (\ref{11}), we obtain
\begin{eqnarray}\label{26}
& &R(X,J_{\alpha}X,Z,J_{\alpha}Z)+R(X,J_{\alpha}X,J_{\beta}Z,J_{\gamma}Z)\\
&+&
R(J_{\beta}X,J_{\gamma}X,Z,J_{\alpha}Z)+ 
R(J_{\beta}X,J_{\gamma}X,J_{\beta}Z,J_{\gamma}Z)\nonumber\\
&=&
\frac{1}{n}\left(\rho_{\alpha}(X,J_{\alpha}X)+ 
\rho_{\alpha}(J_{\beta}X,J_{\gamma}X)\right) 
g(Z,Z) = \frac{2}{n}\rho_{\alpha}(X,J_{\alpha}X)g(Z,Z),\nonumber
\end{eqnarray}
where the last equality of (\ref{26}) is a consequence of the following identity
$$
\rho_{\alpha}(J_{\beta}X,J_{\gamma}X)=-\rho_{\beta}(J_{\beta}X,X)= 
\rho_{\alpha}(X,J_{\alpha}X).
$$
The left side of (\ref{26}) is symmetric with respect to the vectors $X,Z$ because
$D=0$. Hence,
$\rho_{\alpha}(X,J_{\alpha}X)g(Z,Z)=\rho_{\alpha}(Z,J_{\alpha}Z)g(X,X), 
\alpha=1,2,3$.
The last equality together with (\ref{a1}) implies  (\ref{27}). \hfill {\bf Q.E.D.}
\begin{th}\label{th2}
Let $(M,g,(J_{\alpha}))$ be a 4n-dimensional ($n>1$) QKT manifold with 4-form $dT$ 
of type (2,2) with respect to each $J_{\alpha}, \alpha =1,2,3$. Suppose that the 
torsion 
is parallel with respect to the QKT-connection. Then  $(M,g,(J_{\alpha}))$ 
is either a HKT
manifold with parallel torsion or a QK manifold.
\end{th} 
{\it Proof.} We apply Proposition~\ref{teh}. 
If the function $\lambda =0$ then $\rho_{\alpha}=0, 
\alpha =1,2,3$, by (\ref{27}) and Proposition~\ref{p2} implies that the QKT 
manifold is actually a HKT manifold.

 Let $\lambda \not= 0$. The condition (\ref{27}) determines the torsion completely. We proceed
involving (\ref{12}) into the computations as in \cite{HOP}. We calculate using
({\ref{1}) and (\ref{27}) that 
\begin{equation}\label{28}
(\nabla_Z\rho_{\alpha})(X,Y) =\lambda\left\{\omega_{\beta}(Z)F_{\gamma}(X,Y) -
\omega_{\gamma}(Z)F_{\beta}(X,Y)\right\} - d\lambda(Z)F_{\alpha}(X,Y).
\end{equation}
Applying the operator $d$ to (\ref{11}), we get taking into account (\ref{27})
 that  
\begin{equation}\label{29}
d\rho_{\alpha} = \lambda(F_{\beta}\wedge \omega_{\gamma} - \omega_{\beta} \wedge
F_{\gamma})
\end {equation}
On the other hand, we have
\begin{equation}\label{30}
d\rho_{\alpha} =
{\sigma \atop XYZ} \left\{(\nabla_Z\rho_{\alpha})(X,Y)
+\lambda(T(X,Y,J_{\alpha}Z)\right\}, \quad \alpha=1,2,3.
\end{equation}
Comparing the left-hand sides of (\ref{29}) and (\ref{30}) and using (\ref{28}), we
derive
$$
\lambda {\sigma \atop XYZ}\left\{(T(X,Y),J_{\alpha}Z)\right\} = d\lambda\wedge
F_{\alpha}(X,Y,Z), \quad \alpha=1,2,3.
$$
The last equality implies 
$
\lambda T=J_{\alpha}d\lambda \wedge F_{\alpha}, \quad \alpha=1,2,3.
$ 
If $\lambda$ is a non zero constant then $T=0$ and we recover the result of  
 \cite{HOP}. If $\lambda$ is not a constant then there exists a point $p\in
M$ and a neighbourhood $V_p$ of $p$ such that $\lambda \Big |_{V_p} \not=0$.
Then 
\begin{equation}\label{c10}
T=J_{\alpha}d\ln\lambda \wedge F_{\alpha}, \quad \alpha=1,2,3. 
\end{equation}
We take  the trace in (\ref{c10}) to obtain
\begin{equation}\label{c11}
4(n-1)J_{\alpha}d\ln\lambda = 0, \quad \alpha =1,2,3.
\end{equation}
The equation (\ref{c11}) forces $d\lambda =0$ since $n>1$ and consequently $T=0$  
by (\ref{c10}).  Hence, the QKT space is a QK manifold which completes the proof.
\hfill {\bf Q.E.D.}  

On a locally homogeneous QKT manifold the torsion and curvature are parallel and 
Theorem~\ref{th2} leads to the following 
\begin{th}\label{q1}
A (locally) homogeneous 4n-dimensional $(n>1)$ QKT manifold with torsion 4-form $dT$ of type
(2,2) is either (locally) homogeneous HKT space or a (locally) symmetric QK space.
\end{th}
Theorem~\ref{q1} shows that there are no homogeneous (proper) QKT manifolds with 
torsion 4-form of type (2,2) in dimensions greater than four  
which is proved in \cite{homo} by different 
methods using the Lie algebra arguments.

\section{Four dimensional QKT manifolds}

In dimension 4 the situation is completely different 
from that described in Theorem~\ref{th1} and Theorem~\ref{th2} in higher dimensions. 
For a given quaternionic 
 structure on a 4-dimensional manifold $(M,g(H))$ (or equivalently, given 
an orientation and a conformal class of Riemannian metrics \cite{GT}) there are many 
QKT structures \cite{HOP}. More precisely, all QKT structures associated with $(g,(H))$ 
depend on a  1-form $\psi$ due to the general identity 
\begin{equation}\label{tri1}
*\psi=-J\psi\wedge F,
\end{equation} 
where $*$ is the Hodge $*$-operator, $J$ is an $g$-orthogonal almost complex structure 
with K\"ahler form $F$ (see \cite{GT}). Indeed, for any given 1-form $\psi$ we may 
define a QKT-connection $\nabla $ as follows: $\nabla = \nabla^g+\frac{1}{2}*\psi$. 
Conversely, any 3-form $T$ can be represented by $T=-*(*T)$ and the connection given 
above is a quaternionic connection with torsion $T=*\psi$. Hence, a QKT structure on 
a 4-dimensional oriented manifold is a pair $(g,t)$ of a Riemannian metric $g$ 
and an 1-form $t$. The choice of $g$ generates three almost complex structures 
 $(J_{\alpha}),\alpha=1,2,3$, satisfying the quaternionic identities \cite{GT}. The
torsion 3-form T is given by
\begin{equation}\label{v1}
T= *t = t_{\alpha}\wedge F_{\alpha} = t_{\beta}\wedge F_{\beta} = 
t_{\gamma}\wedge F_{\gamma}.
\end{equation}
As consequence of (\ref{tri1}), we obtain $*dT=*d*t =-\delta t$. The last identity 
means that the  torsion 3-form $T$ is closed if and 
only if the 1-form $t$ is co-closed. Thus, in dimension 4 there are many 
strong QKT structures.

In higher dimensions the conformal change of the metric induces a unique QKT 
structure by Theorem~\ref{th3}. We may define a QKT connection corresponding to 
a conformally equivalent metric $\bar{g}=fg$ in dimension 4 by (\ref{z1}) and call 
this conformal QKT transformation. In the compact case, taking the Gauduchon metric
of Theorem~\ref{th3}, we obtain 
\begin{pro}
Let $(M,g,(H),\nabla)$ be a compact 4-dimensional QKT manifold. In the conformal class 
[g] there exists a unique (up to homotety) strong QKT structure conformally 
equivalent to the given one.
\end{pro}
Further, we consider QKT structures with parallel torsion. 
We have
\begin{th}\label{teh2}
A 4-dimensional QKT manifold $M$  with parallel torsion 3-form is 
 a strong QKT manifold,  
the torsion 1-form is parallel with respect to the Levi-Civita connection and $M$ is
locally isometric to the product $N^3\times{\cal R}$, where $N^3$ is a three 
dimensional Riemannian manifold admitting a Riemannian connection $\nabla$ with 
totally skew-symmetric torsion, parallel with respect to $\nabla$.
\end{th}
{\it Proof.} The proof is based on  the following
\begin{lem}\label{tri2}
A 4-dimensional QKT manifold has parallel torsion 3-form if and only if it 
has parallel torsion 1-form with respect to the Levi-Civita connection.
\end{lem}
{\it Proof of Lemma \ref{tri2}.} We calculate using (\ref{v1}) and (\ref{1}) that
\begin{eqnarray}\label{v2}
(\nabla_ZT)(X,Y,U) &=& t_{\alpha}(U)\left(\omega_{\beta}(Z)F_{\gamma}(Y,X) -
\omega_{\gamma}(Z)F_{\beta}(Y,X)\right)\\
 &-&
 t_{\alpha}(X)\left(\omega_{\beta}(Z)F_{\gamma}(Y,U)
- \omega_{\gamma}(Z)F_{\beta}(Y,U)\right)\nonumber\\
&+&
t_{\alpha}(Y)\left(\omega_{\beta}(Z)F_{\gamma}(X,U) -
\omega_{\gamma}(Z)F_{\beta}(X,U)\right)\nonumber\\
&+& 
 F_{\alpha}(Y,U)(\nabla_Zt_{\alpha})X +  F_{\alpha}(X,Y)(\nabla_Zt_{\alpha})U
 -F_{\alpha}(X,U)(\nabla_Zt_{\alpha})Y.\nonumber
\end{eqnarray}
Taking the trace in (\ref{v2}), we obtain
\begin{equation}\label{v3}
\sum_{i=1}^{4}(\nabla_ZT)(X,e_i,J_{\alpha}e_i) =-2(\nabla_Zt_{\alpha})X -
2\left(\omega_{\beta}(Z)t_{\gamma}(X) - \omega_{\gamma}(Z)t_{\beta}(X)\right).
\end{equation} 
Using (\ref{1}), we get
\begin{equation}\label{tre3}
(\nabla_Zt_{\alpha})X = (\nabla_Zt)J_{\alpha}X- 
\left(\omega_{\beta}(Z)t_{\gamma}(X) - \omega_{\gamma}(Z)t_{\beta}(X)\right).
\end{equation}
The equation (\ref{tre3}) and (\ref{v3})  yield
\begin{equation}\label{ser2}
\sum_{i=1}^{4}(\nabla_ZT)(J_{\alpha}X,e_i,J_{\alpha}e_i) = 
2(\nabla_Zt)X, \quad \alpha=1,2,3. 
\end{equation}
Then $\nabla t=0$ since the torsion is parallel. But  $\nabla^gt=\nabla t$ by 
(\ref{5'}) and (\ref{v1}). Hence, $\nabla^gt=0$.

For the converse, we insert (\ref{tre3}) into (\ref{v2}) to get
\begin{equation}\label{ser1}
(\nabla_ZT)(X,Y,U)= F_{\alpha}(Y,U)(\nabla_Zt)J_{\alpha}X +
 F_{\alpha}(X,Y)(\nabla_Zt)J_{\alpha}U +
  F_{\alpha}(U,X)(\nabla_Zt)J_{\alpha}Y,
\end{equation}
since the dimension is equal to four.  If $\nabla^g t=0$ then $\nabla t=0$ 
and (\ref{ser1}) leads to $\nabla T=0$  which proves the lemma. 
\hfill {\bf Q.E.D.}

Lemma~\ref{tri2} shows that $(M,g)$ is locally 
isometric
to the Riemannian product ${\mathcal R}\times N^3$ of a real line and a 
3-dimensional manifold $N^3$ (see e.g. \cite{KN}). Using (\ref{v1}) we see that 
$T(t^{\#},X^{\bot},Y^{\bot})=0$ for every vector fields $X^{\bot},Y^{\bot}$ 
orthonormal to the  vector field $t^{\#}$ dual to the torsion 1-form $t$. 
Hence, the torsion $T$ and therefore the connection $\nabla$ descend to $N^3$. 

In particular, $\delta t=0$ and therefore the QKT structure is strong.
\hfill {\bf Q.E.D.}

As a consequence of Theorem~\ref{teh2}, we recover the following two results 
proved in 
\cite{KV} in the setting of naturally reductive homogeneous 4-manifolds
\begin{th}\label{th4}
A (locally) homogeneous 4-dimensional QKT manifold is locally isometric to 
the Riemannian product ${\mathcal R}\times N^3$ of a real line and 
a naturally reductive homogeneous 3-manifold $N^3$.
\end{th}
\begin{th}\label{tre4}
Let $(M,g)$ be a 4-dimensional compact homogeneous QKT manifold. Then the 
universal covering space $\tilde M$ of $M$ is isometric to the Riemannian 
product 
${\mathcal R}\times N^3$ of a real line and the three dimensional space $N^3$ is 
one of the following

i) $R^3, S^3, {\mathcal H}^3$;

ii) isometric to one of the following Lie groups with a suitable left invariant
metric:

1. $SU(2)$;

2. $\tilde {SL(2,{\mathcal R})}$, the universal covering of $SL(2,{\mathcal R})$;

3. the Heisenberg group.
\end{th}
Theorem~\ref{tre4} is based on the classification of 3-dimensional simply 
connected naturally reductive homogeneous spaces given in \cite{TV}.  

\subsection{Einstein-like QKT 4-manifolds}

It is well known \cite{Ber,A1} that a 4n-dimensional ($n>1)$ QK manifold is Einstein
and the Ricci forms satisfy $\rho_{\alpha}(X,J_{\alpha}Y)= \rho_{\beta}(X,J_{\beta}Y)=
\rho_{\gamma}(X,J_{\gamma}Y)= \lambda g(X,Y)$, where $\lambda$ is a constant.
However, the assumptions that these properties hold on a QKT manifold ($n>1)$ force
the torsion to be zero \cite{HOP} and the QKT manifold is a QK manifold. Actually,
we have already generalised this result proving that if $\lambda $ is  not a 
constant 
 the torsion has to be zero (see the proof of Theorem~\ref{th2}).
 
If the dimension is equal to 4 the situation is different. In this section we show that 
there exists a 4-dimensional (proper) QKT manifold satisfying similar curvature properties
as those mentioned above.

We denote by $K$ the following (0,2) tensor 

$K(X,Y):=
\rho_{\alpha}(X,J_{\alpha}Y)+ \rho_{\beta}(X,J_{\beta}Y)+
\rho_{\gamma}(X,J_{\gamma}Y)$. 

The tensor $K$ is independent of the chosen local
almost complex structures $(J_{\alpha})$ because of the following
\begin{pro}\label{qk1}
Let $(M,g,(J_{\alpha}),\nabla)$ be a 4-dimensional QKT manifold. Then:
\begin{equation}\label{5.67}
K=-Ric +\nabla^g t -\frac{\delta t}{2}g;
\end{equation}
\begin{equation}\label{5.68}
Skew(Ric)= -\frac{1}{4}<dt,F_{\alpha}>F_{\alpha} 
+\frac{1}{2}(d^ct)_{\alpha}, \quad \alpha=1,2,3;
\end{equation}
\begin{equation}\label{5.69}
Ric^g= Sym(Ric) +\frac{1}{2}(|t|^2g -t\otimes t),
\end{equation}
where $<,>$ is the scalar product of tensors induced by $g$, $Skew$
(resp. $Sym$) denotes the skew-symmetric (resp. symmetric) part of a tensor.

In particular, the Ricci tensor is symmetric if and only if the torsion 1-form 
is closed.
\end{pro}
{\it Proof.} We use (\ref{20}). From (\ref{ser2}) and (\ref{tir2}), we obtain
\begin{equation}\label{ap4}
D_{\alpha}(X,J_{\alpha}Y) = (\nabla _Xt)Y - 
(\nabla _{J_{\alpha}Y}t)J_{\alpha}X, \quad \alpha=1,2,3.
\end{equation}
To compute  $B_{\alpha}$ we need the following general identity
\begin{lem}\label{ei3}
On a 4-dimensional QKT manifold we have
$
 {\sigma \atop XYZ}g(T(X,Y),T(Z,U))=0.
$
\end{lem}
{\it Proof of Lemma~\ref{ei3}.} Since ${\sigma \atop XYZ}g(T(X,Y),T(Z,U))$ is a 
4-form it is sufficient to check the equality for a basis of type 
$\{X,J_{\alpha}X,J_{\beta}X,J_{\gamma}X\}$. The last claim is obvious because of (\ref{v1}).

For each $\alpha \in \{1,2,3\}$, 
Lemma~\ref{ei3}, (\ref{ser1}) and (\ref{ser2}) yield
\begin{equation}\label{ap5}
B_{\alpha}(X,J_{\alpha}Y) = 
\sum_{i=1}^4{\sigma \atop XJ_{\alpha}Ye_i}(\nabla_XT)(J_{\alpha}Y,e_i,J_{\alpha}e_i)=
(\nabla_Xt)Y + (\nabla_{J_{\alpha}Y}t)J_{\alpha} X -\delta t g(X,Y).
\end{equation}
Substituting (\ref{ap4}), (\ref{ap5}) into (\ref{20}) and putting $n=1$, we 
derive (\ref{5.67}) since $\nabla^gt=\nabla t$. 
Taking the trace in (\ref{ser1}), we get
$
\sum_{i=1}^4(\nabla _{e_i}T)(e_i,X,Y)= \frac{1}{2}\sum_{i=1}^4 dt(e_i,J_{\alpha}e_i)
F_{\alpha}(X,Y) + dt(J_{\alpha}X,J_{\alpha}Y), \alpha=1,2,3.
$
Then (\ref{5.68}) follows from the last equality and Remark~3. The equation 
(\ref{5.69}) is a direct consequence of (\ref{15}) and (\ref{v1}). \hfill {\bf Q.E.D.}

A $4n$-dimensional
 QKT manifold $(M,g,(J_{\alpha}),\nabla)$ is said to be a {\it Einstein QKT manifold} 
if the symmetric part $Sym(Ric)$ of the Ricci tensor of $\nabla$ is a scalar 
multiple of the metric $g$ i.e.
$
Sym(Ric) =\frac{Scal}{4n}g,
$
where $Scal=tr_g Ric$ is the scalar curvature of $\nabla$.

We note that the scalar curvature $Scal$ of an Einstein QKT manifold may not be 
a constant.

We shall say that a 4-dimensional QKT manifold is {\it sp(1)-Einstein} if 
 the symmetric part  $Sym(K)$ of the tensor $K$ is a scalar 
multiple of the metric $g$ since the tensor $K$ is determined by the sp(1)-part 
of the curvature.
On a sp(1)-Einstein QKT manifold 
$
Sym(K) =\frac{Scal^K}{4}g,
$
where $Scal^K=tr_g K$. 

For a given QKT manifold with torsion 1-form $t$ we consider the corresponding 
Weyl structure $\nabla^W$, i.e. the unique torsion-free linear connection 
determined by the condition 
\begin{equation}\label{qw}
\nabla^Wg=-t\otimes g.
\end{equation} 
Conversely, in dimension 4, 
to a given Weyl structure $\nabla^Wg=\psi\otimes g$ we associate the QKT connection
with torsion $T=*(-\psi)$. 
Note that a given Weyl structure on a conformal manifold (M,[g]) does not depend on 
the particularly chosen metric $g\in [g]$ but 
depends on the conformal class [g]. 
A Weyl structure is said to be {\it Einstein-Weyl} if the 
symmetric part $Sym(Ric^W)$ of its Ricci tensor is a scalar multiple of the metric $g$.
Weyl structures and especially Einstein-Weyl structures have been much studied. 
For a nice overview of Einstein-Weyl geometry see \cite{CP}. The next theorem
shows the link between Einstein-Weyl geometry and sp(1)-Einstein QKT 
manifolds in dimension 4.
\begin{th}\label{qkw}
Let $(M,g,(J_{\alpha}),\nabla)$ be a 4-dimensional QKT manifold with torsion 1-form $t$.
 The following 
conditions are equivalent:

i) $(M,g,(J_{\alpha}),\nabla)$ is a sp(1)-Einstein QKT manifold. 

ii) The corresponding Weyl structure is an Einstein-Weyl structure.
\end{th}
{\it Proof.} The Weyl connection $\nabla^W$ determined by (\ref{qw})
is given explicitly by
$$
\nabla ^W_X Y= \nabla^g _X Y + \frac{1}{2} t(X)Y
+ \frac{1}{2}t(Y)X - \frac{1}{2} g(X,Y)t^{\#}.
$$

The symmetric part of its Ricci tensor is equal to
\begin{equation}\label{wzl1}
Sym(Ric^W) = Ric^g - Sym(\nabla^g t) -\frac{1}{2}(|t|^2g -t\otimes t) +
\frac{\delta t}{2}g.
\end{equation}
Keeping in mind that $\nabla^g t=\nabla t$, we get from (\ref{5.67}), (\ref{5.69}) 
and (\ref{wzl1}) that \\
 $Sym(Ric^W) =- Sym(K)$. The theorem follows 
from the last equality. \hfill {\bf Q.E.D.}

It is well known \cite{AHS,Sal4} 
that on a 4-dimensional conformal manifold there exists a hypercomplex structure iff 
the conformal structure has anti-self-dual Weyl tensor (see also \cite{GT}). 
Every 4-dimensional hypercomplex manifold $(M,g,(H_{\alpha}))$, i.e. 
(an oriented anti-self-dual 4-manifold) carries a unique HKT structure in view of the
results in \cite{GT,Ga1}. Indeed, let $\theta =\theta_{\alpha}=\theta_{\beta}=
\theta_{\gamma}$ be the common Lee form. The unique HKT structure is defined by 
$\nabla =\nabla^g-\frac{1}{2}*\theta$ \cite{GT} (the uniqueness is a consequence of 
a general result in \cite{Ga1}, see also \cite{GP}). 
The HKT structure on a 4-dimensional hypercomplex manifold is sp(1)-Einstein 
since the tensor $K$ vanishes.
 The corresponding Weyl structure to the given HKT structure on a 
4-dimensional hyperhermitian manifold is the Obata connection \cite{GT}, i.e. the 
unique torsion-free linear connection which preserves each of the three hypercomplex 
structures. As a consequence of Theorem~\ref{qkw}, we recover the result 
in \cite{PS} which states that the Obata connection of a hyper complex 4-manifold
is Einstein-Weyl and the symmetric part of its Ricci tensor is zero.

Theorem~\ref{qkw} and (\ref{5.67}) show that every Einstein-Weyl structure 
determined by (\ref{qw})
on a 4-dimensional 
conformal manifold whose vector field dual to the  1-form $t$ is Killing, induces 
an Einstein and sp(1)-Einstein QKT structure. 
\begin{co}
Let $(M,[g],\nabla^W)$ be a compact 4-dimensional Einstein-Weyl manifold. Then the 
corresponding QKT structure to the Gauduchon metric of $\nabla^W$ is Einstein and 
sp(1)-Einstein.
\end{co}
{\it Proof.} On a compact Einstein-Weyl manifold the vector field dual to the Lee 
form of the Gauduchon metric is Killing by the result of Tod \cite{tod}. Then the 
claim follows from Theorem~\ref{qkw} and (\ref{5.67}). \hfill {\bf Q.E.D.}

The Ricci tensor of a 4-dimensional QKT manifold is symmetric iff the torsion 1-form is
closed by Proposition~\ref{qk1}. Applying Theorem~3 in \cite{G2} and using 
Theorem~\ref{qkw}, we obtain
\begin{co}\label{qwe}
Let $(M,g,(J_{\alpha}),\nabla)$ be a 4-dimensional compact  sp(1)-Einstein
QKT manifold with symmetric Ricci tensor. Suppose that the torsion 1-form is not exact.
Then the torsion 1-form corresponding to the Gauduchon metric $g_G$ of 
$(M,g,(J_{\alpha}),\nabla)$ is parallel with respect to the Levi-Civita connection 
of $g_G$ and the universal cover of $(M,g_G)$ is isometric to 
${\mathcal R}\times S^3$. 
 In particular, the quaternionic bundle $(J_{\alpha})$ admits hypercomplex structure.
\end{co}

A lot is known about Einstein-Weyl manifolds (see a nice survey \cite{CP}). There
are many (compact) Einstein-Weyl 4-manifolds (e.g. $S^2\otimes S^2$). Among
them there are (anti)-self-dual as well as non (anti)-self-dual. We mention here
the Einstein-Weyl examples of Bianchi IX type metric \cite{Bi,Ma,MaP}. All these 
Einstein-Weyl 4-manifolds admit sp(1)-Einstein QKT structures by Theorem~\ref{qkw}.

It is also known that there are obstructions to the existence of Einstein-Weyl 
structures 
on compact 4-manifold \cite{PPS1}. If the manifold $M$ is finitely covered by 
$T^2\otimes S^2$ which cannot be Einstein-Weyl then $M$ does not admit
Einstein-Weyl structure and  therefore there are no sp(1)-Einstein
structures on $M$.

\bigskip {\bf Authors' address:}\\[2mm]
Stefan Ivanov\\ University of Sofia, Faculty of
Mathematics and Informatics, Department of Geometry, \\ 5 James
Bourchier blvd, 1126 Sofia, BULGARIA.\\ 
E-mail: {\tt ivanovsp@fmi.uni-sofia.bg}
\end{document}